\documentclass[a4paper,11pt]{amsart}
\usepackage[all]{xy}
\usepackage{amsmath,amsthm,amsfonts,amssymb,color,graphicx,
hyperref}

\newenvironment{nouppercase}{%
  \renewcommand{\uppercasenonmath}[1]{}}{}

\begin{document}
\thispagestyle{empty}

\title[\sc Degree of the exceptional component of foliations]
  {Degree of the exceptional component of foliations\\ of degree two and codimension one in  \mbox{$\mathbb P^ 3$}}
  
\author[A. Rossini]{Artur Rossini}
\author[I. Vainsencher]{Israel Vainsencher}


\begin{nouppercase} \maketitle \end{nouppercase}
\begin{abstract}
The purpose of this work is to obtain the
degree of the exceptional component of the space
of holomorphic foliations of degree two and codimension
one in $\mathbb{P}^3$. This component is the
closure of the orbit 
of the foliation defined by the differential form
$$\omega=(3fdg-2gdf)\big/x_0,
\ \mbox{ where } \ f=x_0^2x_3-x_0x_1x_2+\dfrac{x_1^3}{3},\ 
g=x_0x_2-\dfrac{x_1^2}{2}
$$
under the natural action of the group of automorphisms of $\mathbb{P}^3$.
Our first task is to unravel a geometric characterization
of the pair $g,f$. This leads us to 
the construction of a parameter space as 
an explicit fiber bundle over the 
variety of complete flags.
Using tools from equivariant intersection theory,
especially Bott's formula,
the  degree is expressed as an integral over our
parameter space.
\end{abstract}

\section{Introduction}

A holomorphic foliation of codimension one and degree $d$
in the complex projective space $\mathbb{P}^n$ is given
by a differential 1-form 
$\omega = A_0dx_0 + \cdots + A_ndx_n$
where $A_0,\dots,A_n$ are homogeneous polynomials of degree $d+1$, satisfying the conditions (i) $A_0x_0 + \cdots +A_nx_n = 0$ and (ii) (integrability) $\omega \wedge d\omega = 0.$

These two conditions define a closed subscheme
$\overline{\mathcal{F}(d,n)}$ of the projective
space of global sections of the twisted cotangent sheaf
$\Omega_{\mathbb{P}^n}^1(d+2)$. 
Naturally, $\overline{\mathcal{F}(d,n)}$ possesses a
decomposition into irreducible components.  This
decomposition is completely known in cases $d = 0$ and $d
= 1$, thanks to the pioneering book by
\cite{Jouanolou}. For $n \geq 3$, $d=2$ 
a full description can be found in the celebrated work by
\cite{Alcides} where we learn that there are six
components of $\overline{\mathcal{F}(2,3)}$,
among which
is the \textit{exceptional component} $E(3)$,
our main object of interest here. 

There are a few other known components for $d \geq 3$,
such as the pullback, the rational  and
the logarithmic components, see \cite{NetoeEndonvein},
\cite{CukiermanVI} and \cite{Omegar}. In
\cite{CukiermanVI} the authors managed to compute the
degree of some rational components. In particular, the
degrees of the rational components $\mathcal R(n,2,3)$
were computed for $n \leq 5$; in \cite{daniel} this
is extended to the
rational components $\mathcal R(n,2,2r+1)$ for arbitrary $n,r$. 

The \textit{exceptional component} $E(3)$ is the
closure of the orbit, under the natural action of
Aut($\mathbb{P}^3$), of the foliation defined by the differential form
\begin{equation}\label{fgi}
\omega_0=\frac{3f_0dg_0-2g_0df_0}{x_0},
\ \mbox{ where } \ f_0=x_0^2x_3-x_0x_1x_2+\dfrac{x_1^3}{3},\ 
g_0=x_0x_2-\dfrac{x_1^2}{2}
.\end{equation}

Our goal is to describe the geometry of the family
of pairs $g,f$ in the orbit closure, enabling us to
compute the degree of $E(3)$. 
Our strategy  starts with a
complete flag in $\mathbb{P}^3$ 
$$p \textrm{ (point)} \in \ell \textrm{ (line) } 
\subset v \textrm { (plane)},
$$
over which we describe suitable cubic forms $f$ and 
quadratic forms $g$ that fulfill some special
conditions. The polynomials $f,g$ will 
give us the differential 1-form
$$\omega = \dfrac{3fdg-2gdf}{h},
$$
where $h$ is an equation of the plane $v$ of the
given flag. 

The main technical difficulty is 
to solve  the indeterminacies of the rational map
$$(g,f) \mapsto \omega.
$$ This is done by a careful analysis of the irreducible
components of the  indeterminacy locus, blown up one at a
time. Macaulay2 (\cite{M2}) is extensively used.
We get a description of a parameter
space for the exceptional component, over which we can apply appropriate tools from equivariant intersection theory to compute the desired degree.

\section{A parameter space 
	for the Exceptional Component}\label{description}

The  component $E(3)$ 
is the orbit closure of the foliation defined by the
1-form $\omega_0$ as in (\ref{fgi}).
The dimension of  $E(3)$
is equal to 13, see  \cite{jinvariante}.

Explicitly, we have
\begin{equation}\label{omega0}
\begin{array}{rcrcr}
\omega_0&=&
(x_1x_2^2-2x_1^2x_3+x_0x_2x_3)dx_0& +&
x_0(3x_1x_3-2x_2^2)dx_1  \\
&&\! +x_0(x_1x_2-3x_0x_3)dx_2& +&\!\!
x_0(2x_0x_2-x_1^2)dx_3.
\end{array}
\end{equation}
The singular locus of the foliation (\ref{omega0}) consists
of a union of three curves:

\begin{enumerate}
	\item the conic given by ideal     
	$\langle x_0,x_2^2-2x_1x_3 \rangle$
	(it lies in  the plane $x_0=0$);
	
	\item the line $\ell_0$ defined by
	$ x_0=x_1=0$;
	
	\item the twisted cubic given by $ \langle 2x_2^2-3x_1x_3, \ x_1x_2-3x_0x_3, \ x_1^2-2x_0x_2 \rangle$. 
\end{enumerate} \medskip

These three components meet at the
point $p_0:=(0:0:0:1) \in \mathbb{P}^3$. 
Let us examine the geometry of the  surface defined
by the cubic form
$f_0$.
We see  that it is an irreducible cubic, singular along the
line \,$\ell_0$.
Moreover, given a point $ p_t = (0: 0: t: 1) \in \ell_0
$, the tangent cone to this cubic surface at the point $
p_t $ has equation
$$ x_0^2 - tx_0x_1 = x_0 (x_0-tx_1). 
$$
Thus,  the tangent cone is a pair of planes containing the double
line \,$\ell_0$, one of which is fixed and the other varies
with the point $ p_t$.

We have also perceived that at  the special
point $ p_0= (0: 0: 0: 1) $ the tangent
cone to the cubic is the double plane $ x_0^2 = 0$.

Therefore, the cubic $f_0$ comes endowed with
a companion complete flag 
\begin{equation} \label{pdv}
\varphi_0:\ 
p_0=\{x_0=x_1=x_2=0\} \in \ell_0 = \{x_0=x_1=0\} 
\subset v_0 =
\{x_0=0\} .
\end{equation} 

As for the quadric
$$ g_0 = x_0x_2 - \frac{x_1^2}{2},$$
we have a cone containing the line $\ell_0$
and vertex 
$p_0$. Moreover, the tangent plane to the cone at a
smooth point on $\ell_0$ is precisely the same plane $v_0$.

\subsection{\bf Remark.}	\label{recupflag}  
	The flag can move under the action of an
	automorphism of $\mathbb{P}^3$, but it can be
	recovered directly from the 1-form $\omega $
	defining the exceptional foliation, just by
	looking at the singular locus. For any such
	$\omega$, the singular locus has three components
	- a conic living on the new plane, the new line
	and a twisted cubic that meets the other two
	components just at the new point; the line is
	tangent to the twisted cubic, the plane is the
	osculating plane and the conic is the osculating
	conic. By the way, the dimension 13 mentioned just
	above is the dimension of the family of pointed
	twisted cubics.\medskip

These considerations about the cubic/quadric pairs
as in (\ref{fgi}) lead us to consider the construction of a
parameter space for the family of exceptional foliations
as a fiber bundle over 
the variety $ {\mathbb F}$  of complete flags on $ \mathbb{P}^3 $, 
$$
p \textrm{(point)} \in \ell \textrm{(line)} \subset \pi
\textrm{(plane)}.
$$
Given such a flag we take: \medskip

$\bullet$ A cubic surface $f$ with the properties:
\begin{enumerate}
	\item $f$ is singular along $\ell$;   
	\item the tangent cone to $f$ at
	each point $q \in \ell$ is
	the union of two planes $\pi \cup \alpha_q$.  
	\item $\alpha_p = \pi$, that is, at the point $p$
	the tangent cone  to $f$ is the double plane.
\end{enumerate} 

$\bullet$ A quadratic cone $g$ with the properties:
\begin{enumerate}
	\item the line $\ell$ is contained in the cone $g$;
	\item the plane $\pi$ is tangent to $g$ along
	the line $\ell$;
	\item the point $p$ is in the vertex of
	the cone.  
\end{enumerate}

An exceptional foliation is then given by the differential form
$$\omega = \dfrac{3fdg - 2gdf}{h},$$
where $h$ is an equation of the plane $\pi$. 

\subsection{\bf Definition.}
	For simplicity we will call here a cubic or a quadric that meet the requirements above as \textit{special}.\medskip

We realize that after the
construction of  special pairs
$(f, g)$, it is still necessary to
impose the condition that the differential form $
3fdg-2gdf $ be divisible by the equation of the plane. We
refer to this as  the \textit{divisibility condition}.
In order to better understand this divisibility
condition, let us fix a complete flag as in
(\ref{pdv}). 

The requirements about special cubics/quadrics show that
we have
\begin{equation}\label{ff}\left\{\begin{array}l
f=\left(a_0x_0 + a_1x_1 + a_2x_2+a_3x_3\right)x_0^2 +
a_4x_0x_1^2+a_5x_0x_1x_2 + a_6x_1^3\\
g=b_0x_0^2 + b_1x_0x_1 + b_2x_0x_2 + b_3x_1^2.
\end{array}\right.
\end{equation}

We register the following invariant description.

\subsection{\bf Lemma.} \label{monsfg} 
	The $7$ monomials appearing in \ $f$ \medskip
	\\\centerline{$x_{0}^{3},x_{0}^{2}x_{1},
		x_{0}^{2}x_{2}, x_{0}^{2}x_{3}, x_{0}x_{1}^{2},   x_{0}x_{1}x_{2}, x_{1}^{3} $}
	are the monomial generators of the
	subspace  of cubics,
	\\\centerline{$x_0 ^2\langle x_0,\dots,x_3\rangle+x_0
		\langle x_0,x_1\rangle \langle x_0,x_1,x_2\rangle + \langle x_0,x_1\rangle ^3\subset S_3$.}
	Likewise, the $4$ monomials \ $x_0^2,x_0x_1,x_0x_2,x_1^2$ \ 
	in the quadric \ $g$ \ generate the subspace
	\\\centerline{$x_0\langle x_0,x_1,x_2\rangle +  \langle x_0,x_1\rangle ^2\subset S_2
		$.}\medskip

Varying the flag, we obtain equivariant vector
subbundles 
\begin{equation}
\label{AB}
\mathcal A\subset S_3 \ \text{ and } \ \mathcal B\subset S_2.
\end{equation}
with respective ranks 7 and 4.
Projectively we obtain, at this stage, the variety
$$
\mathbb{P}(\mathcal A) \times_{\mathbb F}\mathbb{P}(\mathcal B),
$$
a $\mathbb{P}^6 \times \mathbb{P}^3$ bundle
of pairs $(f,g)$ over a fixed flag. 

Continuing the discussion about the divisibility
condition, we may write
$$3fdg-2gdf = x_0\omega_1 + \left[(3a_6b_1-2a_4b_3)x_1^4 + (3a_6b_2 - 2a_5b_3)x_1^3x_2\right]dx_0,$$
where $\omega_1$ is a 1-form. From this it follows that the
divisibility condition is given, on the fiber over our
fixed flag $\varphi$,  by the equations 
$$\left\{
\begin{array}{c}
3a_6b_1 = 2a_4b_3, 
\\3a_6b_2 = 2a_5b_3
\end{array}
\right.$$

This locus consists of two irreducible components inside
$\mathbb{P}^6 \times \mathbb{P}^3$, both of
codimension two:
\begin{equation}
\label{divcond1}
\begin{array}{lrc}
&  (\star)& 
\left\{\begin{array}{lcr}
3a_6b_1 &=& 2a_4b_3\\
3a_6b_2 &=& 2a_5b_3\\
a_4b_2 &=& a_5b_1
\end{array}\right\} \end{array} \ \ \ \text{ and } \ \ \ \ 
\begin{array}{lcr}
& (\star \star )& a_6=b_3=0.
\end{array}
\end{equation}

For a pair $(g,f)$ satisfying $(\star \star )$,
we actually get
$$\left\{\begin{array}{l}
f=x_0\left(a_0x_0^2 + a_1x_0x_1 + a_2x_0x_2+a_3x_0x_3 +
a_4x_1^2+a_5x_1x_2\right),\\
g=x_0\left(b_0x_0 + b_1x_1 + b_2x_2\right).
\end{array}\right.$$
It means that a general element in the second component $(\star \star)$
consists of a cubic and a quadric both
divisible by the equation of the plane,
certainly not interesting for our study of the 
exceptional component.

\medskip
Henceforth, we refer to the {\em divisibility condition}
as  the equations $(\star)$ in (\ref{divcond1}). \medskip

Let $\mathbb{G}=\mathbb{G}(1,3)$ be the Grassmann variety of lines in $\mathbb{P}^3$, with tautological sequence
$$0 \longrightarrow \mathcal{S} \longrightarrow \mathbb{G} \times \mathbb{C}^4 \longrightarrow \mathcal{Q} \longrightarrow 0,$$
where rank $(\mathcal S)$=2.
Denote by  $\mathcal{Q}^{\vee}$ the dual
of $\mathcal{Q}$. We also have the tautological sequence
over $\mathbb{P}^3$,
$$
0 \longrightarrow \mathcal{O}_{\mathbb{P}^3}(-1)
\longrightarrow  \mathbb{P}^3
\times \mathbb{C}^4 \longrightarrow \mathcal{P} \longrightarrow 0
$$
Then,
$$
\begin{array} {*5c}
\mathbb{P}(\mathcal{S})&=&\{(p,\ell) \,|\, p \in \ell
\} &\subset& \mathbb{P}^3 \times \mathbb{G}
\\
\mathbb{P}(\mathcal{Q}^{\vee})&=&
\{(\ell, \pi) \,|\, \ell \subset \pi \}&
\subset& \mathbb{G} \times \check{\vphantom{|}\mathbb{P}}^3
.\end{array}$$
The variety of complete flags  is just the fiber product
$$\mathbb{F}=\mathbb{P}(\mathcal{S}) \times _\mathbb{G}
\mathbb{P}(\mathcal{Q}^{\vee}) = \{(p,\ell,\pi) | \ p \in
\ell \subset \pi \}
\subset\mathbb{P}^3\times \mathbb{G}\times \check{\vphantom{|}\mathbb{P}}^3.
$$
The tautological bundles of these spaces lift
to bundles over $\mathbb{F}$ still denoted by the
same letters. They
fit together into the diagram (pullbacks omitted),
$$\begin{array}{*9c}
\mathcal O_{\mathcal Q^\vee}(-1)&=&
\mathcal O_{\check{\vphantom{|}\mathbb P}^3}(-1) & \hookrightarrow & \mathcal Q^\vee & 
\hookrightarrow & \mathcal P^\vee & 
\hookrightarrow & S_1:=(\mathbb{C}^4)^{\vee}\\ &\varphi_0:&
\langle x_0 \rangle &\subset& \langle x_0,x_1\rangle&\subset&\langle x_0,x_1,x_2\rangle
&\subset&\langle x_0,\dots,x_3 \rangle
\end{array}$$
where the bottom row indicates the corresponding
fibers over the flag (\ref{pdv}). \medskip

In view of Lemma \ref{monsfg}, we have  the surjections
$$\begin{array}{c}
\left(\mathcal O_{\mathcal Q^\vee}(-2)\otimes S_1\right)\oplus\left(
\mathcal O_{\mathcal Q^\vee}(-1)\otimes{\mathcal Q^\vee}\otimes\mathcal
P ^\vee\right)
\oplus{}{\textrm{Sym}}_3\mathcal Q^\vee\twoheadrightarrow
\mathcal A\subset S_3
\\
(\mathcal O_{\mathcal Q^\vee}(-1)\otimes\mathcal P^\vee)
\oplus{}\textrm{Sym}_2\mathcal Q^\vee
\twoheadrightarrow
\mathcal B\subset S_2.
\end{array}
$$

By construction, the vector bundles $\mathcal A$ and $\mathcal B$ fit into the
exact sequences 
$$\begin{array} {*4cl}
S_1\otimes\mathcal{O}_{\ensuremath{\check{\vphantom{|}\mathbb P}}^3}(-2)& \hookrightarrow & \mathcal A&
\twoheadrightarrow&
\overline{\mathcal A}:=\mathcal A\big/\left(S_1\otimes\mathcal{O}_{\ensuremath{\check{\vphantom{|}\mathbb P}}^3}(-2)
\right)\\
\mathcal O_{\ensuremath{\check{\vphantom{|}\mathbb P}}^3}(-2)&\hookrightarrow&
\mathcal B &\twoheadrightarrow &\overline{\mathcal B}:=\mathcal B\big/\left(\mathcal{O}_{\ensuremath{\check{\vphantom{|}\mathbb P}}^3}(-2)
\right).
\end{array}$$

\subsection{\bf Lemma.}\label{ovAovB}
	$\overline{\mathcal A}$ is isomorphic to 
	$\left(\mathcal Q^\vee/\mathcal
	O_{\ensuremath{\check{\vphantom{|}\mathbb P}}^3}(-1)\right)\otimes\overline{\mathcal B}$.

\begin{proof}
	Indeed,
	on the fiber over $\varphi_0$ we have
	$$(
	\mathcal Q^\vee/\mathcal O_{\ensuremath{\check{\vphantom{|}\mathbb P}}^3}(-1))_{\varphi_0}=\langle x_0,x_1 \rangle/ \langle x_0\rangle
	= \langle \overline{x_1}\rangle 
	$$
	hence
	$$
	\left((
	\mathcal Q^\vee/\mathcal O_{\ensuremath{\check{\vphantom{|}\mathbb P}}^3}(-1))\otimes\overline{\mathcal B}\right)_{\varphi_0}=
	\langle \overline{x_1}\rangle \otimes
	\langle \overline{x_0x_1},\overline{x_0x_2},\overline{x_1^2} \rangle
	=\overline{\mathcal A}_{\varphi_0},
	$$
	where the = signs mean isomorphisms of
	representations of the stabilizer of the flag
	$\varphi_0$. 
	\end{proof}

\subsection{\bf Definition.}
	We define
	$\mathbb{X}=\mathbb P(\mathcal B)$
	the corresponding projective bundle of special quadrics.
	The fiber $\mathbb{X}_{\varphi_0}
	$ is the $\mathbb{P}^3$ of special quadrics over the flag
	$\varphi_0 =(p_{0},\ell_0,v_0)$ fixed in (\ref{pdv}). \medskip

The divisibility condition $(\star)$\,(\ref{divcond1}) can be rewritten as a 
system of linear equations 
in variables $(\underline{a})=(a_4:a_5:a_6)$ with
coefficients $(\underline{b})=(b_0:b_1:b_2:b_3)\in \mathbb{P}^3$,
\begin{equation}\label{syslin}
\left[\begin{array}{ccc}
2b_3 & 0     & -3b_1\\
0    & -2b_3 & 3b_2\\
-b_2 & b_1   & 0
\end{array}\right] \cdot 
\left[\begin{array}{c}
a_4\\
a_5\\
a_6
\end{array}\right] =
\left[\begin{array}{c}
0\\
0\\
0
\end{array}\right].
\end{equation}

The matrix of the coefficients has determinant zero and 
generic rank two. This rank drops when
$b_1=b_2=b_3=0\leftrightarrow x_0^2$. Off the
point $x_0^2$  the solution space is spanned by
the vector product of two rows, 
$$(-3b_1b_2:-3b_2^2:-2b_2b_3)=( 3b_1:3b_2: 2b_3)\ .$$

The idea now is to describe the locus of special cubic
forms with the divisibility condition as a
$\mathbb{P}^4$-bundle over $\mathbb{X}$ (since it has
codimension two in the $\mathbb{P}^6$ of special
cubics). But to make it possible, we need to replace
$\mathbb{X}$ by a new space, $\mathbb{X}'$,
for which the rank of the
coefficients matrix in (\ref{syslin}) is two everywhere.

Precisely, think of the fiberwise solution space to
(\ref{syslin}) as defining a rational map 
$\psi:\mathbb X \rightarrow \mathbb{P}(\overline{\mathcal A})$, which on the fiber over the standard
flag $\varphi_0$ (\ref{pdv}) reads
\begin{equation}\label{b2a}
\psi:\ 
(\underline b)\mapsto
(a_4:a_5:a_6)=
(3b_1:3b_2:2b_3).
\end{equation}
Look at the closure $\mathbb{X}'$ of the graph of
$\psi$. On the fiber over $\varphi_0$,
this is the blowup of $\mathbb{P}^3=\mathbb{P}(\mathcal B_{\varphi_0})$ 
at the point $x_0^2$. 
As a matter of fact, this turns out to be the
restriction to the fiber over $ \varphi_0$  
of the blowup of $\mathbb{X}$ along the section
$\mathbb{P}\left(\mathcal{O}_{\ensuremath{\check{\vphantom{|}\mathbb P}}^3}(-2)\right)
\hookrightarrow \mathbb{P}(\mathcal{B})$ over $\mathbb
F$. We have the diagram
$$
\xymatrix
@R2.5pc
{&&
	\,  \mathbb X'
	\ar@{->}[d]_{\text{\small (blowup)}}
	\ar@{->}[drr]^{\ \ \ \ \psi':\text{\small\,\,($\mathbb{P}^1$-bundle)}}
	&&&
	\\ 
	\mathbb{P}\left(\mathcal{O}_{\ensuremath{\check{\vphantom{|}\mathbb P}}^3}(-2)\right)
	\ \ \ar@{=}[drrr] 
	\!\! \ar@{^(-}[r]& 
	\mathbb{P}(\mathcal B) \ar@{=}[r]& \mathbb X
	\ar@{->}[dr]
	\ar@{-->}[rr]^{\!\hskip-.31cm^\psi}&&\mathbb{P}(\overline{\mathcal A})
	\ \ar@{->}[dl] \ar@{=}[r]& \mathbb{P}(\overline{\mathcal B})
	\\&&&\ \ \ \mathbb F \ \ \ & &
}
$$%
where the leftmost inclusion is defined by
taking the square of the equation of the plane in
the flag.
The equality $\mathbb{P}(\overline{\mathcal A})=\mathbb{P}(\overline{\mathcal
	B})$ comes from Lemma\,\ref{ovAovB}; under this
identification,
$$\mathcal O_{\overline{\mathcal A}}(-1)=
\mathcal O_{\overline{\mathcal B}}(-1)\otimes (
\mathcal Q^\vee/\mathcal O_{\ensuremath{\check{\vphantom{|}\mathbb P}}^3}(-1)).
$$
Pulling back the above tautological line
subbundle via $\psi'$ , we get the diagrams
\begin{equation}\label{ovAB}
\begin{array}{*5c}
S_1\otimes\mathcal{O}_{\ensuremath{\check{\vphantom{|}\mathbb P}}^3}(-2)
&\hookrightarrow&\mathcal A'&
\twoheadrightarrow&
\mathcal O_{\overline{\mathcal A}}(-1)\\
||&&\downarrow&&\downarrow\\
S_1\otimes\mathcal{O}_{\ensuremath{\check{\vphantom{|}\mathbb P}}^3}(-2) &\hookrightarrow& {\mathcal A}&
\twoheadrightarrow&
\overline{\mathcal A}\
\end{array}
\end{equation}
and
\begin{equation}\label{ovAB1}
\begin{array}{*5c}
\mathcal{O}_{\ensuremath{\check{\vphantom{|}\mathbb P}}^3}(-2)
&\hookrightarrow&\mathcal B'&
\twoheadrightarrow&
\mathcal O_{\overline{\mathcal B}}(-1)\\
||&&\downarrow&&\downarrow\\
\mathcal{O}_{\ensuremath{\check{\vphantom{|}\mathbb P}}^3}(-2) &\hookrightarrow& {\mathcal B}&
\twoheadrightarrow&
\overline{\mathcal B}\
\end{array}
\end{equation}
where rank $\mathcal A'=5$, rank $\mathcal B'=2$. We have 
\begin{equation}\label{defX'}
\mathbb X'=\mathbb{P}(\mathcal B'),
\end{equation} a $\mathbb{P}^1$ bundle over the $\mathbb{P}^2$ bundle 
$\mathbb{P}(\overline{\mathcal B})$ over the flag variety
$\mathbb F$. 
Define
\begin{equation}\label{Y}
\mathbb Y=\mathbb{P}(\mathcal A'),
\end{equation}
an equivariant $\mathbb{P}^4$--bundle over $\mathbb X'$.

\subsection{\bf Proposition.}
	A general point in $\mathbb Y$ corresponds to a pair
	$(g,f)$ of special quadric, cubic satisfying the
	divisibility condition.
\begin{proof}
	The assertion follows from the previous
	considerations. 
\end{proof}

The next step is to solve the indeterminacies of the rational map
\begin{equation}\label{mapomega}
\begin{array}{ccc}
\mathbb{Y} & \makebox{\Large$\dashrightarrow$}    &
E(3)\subset
\mathbb{P}(H^0(\Omega^1_{\mathbb{P}^3}(4)))\\
(g,g',f) & \longmapsto & \omega = \dfrac{3fdg - 2gdf}{x_0} \ .\\
\end{array}
\end{equation}
This will be accomplished by a sequence of 4 blowups,

\medskip
\medskip
$$\xymatrix
{\mathbb Y_4\,\ar@/^2.0pc/[rrrrr]\,\ar@{->}[r]
	&
	\mathbb Y_3\,\ar@{->}[r]&\mathbb Y_2\,\ar@{->}[r]&\mathbb Y_1\,\ar@{->}[r]&\mathbb Y\,\ar@{-->}[r]&E(3).
}$$

\subsection{\bf Remark.}
	Since all constructions performed so far are equivariant, we drop the reference to the fiber over $\varphi_0$, and simplify notation writing
	$$\mathbb X=\mathbb
	X_{\varphi_0},\ \mathcal A=\mathcal A_{\varphi_0},\dots,\  etc.
	$$

In view of Lemma \ref{ovAovB}, the rational map (\ref{b2a}) can also be written as the rational linear projection map
$$\begin{array}{cccc}
\psi_{\overline{\mathcal{B}}}: & \mathbb{X} &
\dashrightarrow   & \mathbb{P}(\overline{\mathcal{B}})\\
& g:=b_0x_0^2 + b_1x_0x_1 + b_2x_0x_2 + b_3x_1^2 &
\longmapsto
& g':=b_1\overline{x_0x_1} + b_2\overline{x_0x_2} + b_3\overline{x_1^2}
\end{array}$$

(class mod.\,$x_0^2$). Write for short 
\begin{equation} \label{X'}
\mathbb{X}' = \{ (g,g') =
\left((b_0:b_1:b_2:b_3),(u_1:u_2:u_3)\right)
| \ b_iu_j=b_ju_i  \},
\end{equation}
where
\begin{equation} \label{defg'}
g':=u_1\overline{x_0x_1} + u_2\overline{x_0x_2} + u_3\overline{x_1^2}\ .
\end{equation}

Here, the $\underline{u}$ are 
homogeneous coordinates in $\mathbb{P}^2$ = fiber of $\mathbb{P}(
\overline{\mathcal B})$ over $\varphi_0$.
The map $\psi'$ is (see  (\ref{b2a}))
$$\begin{array}{cccc}
\psi': & \mathbb{X'} & \longrightarrow    & \mathbb{P}(\overline{\mathcal{A}})\\
& (g,g')
& \mapsto & 3u_1\overline{x_0x_1^2} + 3u_2\overline{x_0x_1x_2}
+ 2u_3\overline{x_1^3} \\
\end{array}$$

Notice that the divisibility condition
(\ref{divcond1})
has now the expression
\begin{equation}\label{divcond}
\left\{\begin{array}{lcl}
3a_6u_1 &=& 2a_4u_3\\
3a_6u_2 &=& 2a_5u_3\\
a_4u_2 &=& a_5u_1
\end{array}\right.\end{equation}
It can be seen that the system above is of rank 2 for all $(u_1:u_2:u_3)\in\mathbb{P}^2$.
Equations (\ref{divcond}) enables to find all special
cubics over any point $(g,g') \in \mathbb{X}'$.

\section{Solving the Indeterminacies}\label{solvingmap}

Let us work over an affine cover of $\mathbb{X}'$ to deal with the indeterminacies of the map $(g,g',f) \mapsto \omega$. 

Notice that over a fixed flag, a point in $\mathbb{Y}$ has 7 affine coordinates. 

Notation as in (\ref{X'}), consider the open dense set of $\mathbb{X}'$ given by $b_0=u_1=1$. The interesting equations are
$$ \left\{\begin{array}{l}
b_2 = b_1u_2\\
b_3 = b_1u_3
\end{array}\right\} \ \textrm{ and } \ \left\{\begin{array}{l}
3a_6 = 2a_4u_3\\
a_5 = a_4u_2
\end{array}\right\},$$
where the former corresponds to the blowup of our $\mathbb{P}^3$ of
special quadrics at $x_0^2$ (\ref{X'}) and the latter 
to divisibility $(\star)$\,(\ref{divcond}). So,
$$\left\{\begin{array} l
f=\left(a_0x_0 + a_1x_1 + a_2x_2
+a_3x_3\right)x_0^2 +a_4\left(
x_0x_1^2+u_2x_0x_1x_2 + \frac{2}{3}u_3x_1^3\right)
\\
g=x_0^2 + b_1x_0x_1 + b_1u_2x_0x_2 + b_1u_3x_1^2\ .
\end{array}\right.
$$

Computing $\omega$ as in (\ref{mapomega}), we find using Macaulay2
\begin{eqnarray*}
	\omega = [(3a_0b_1-2a_1)x_0^2x_1+(6a_0b_1u_3+a_1b_1-4a_4)x_0x_1^2+ \\ + (3a_0b_1u_2-2a_2)x_0^2x_2-2a_3x_0^2x_3+ \cdots ]dx_0 + \cdots
\end{eqnarray*}
For  $a_0=0$  the 1-form $\omega $ vanishes only if
$a_0=a_1=a_2=a_3=a_4=0,$ which is impossible. Thus, in order to study the indeterminacy locus we can take $a_0=1$.

Setting $a_0=1$ and collecting coefficients of $\omega$, we get a non-reduced and reducible scheme, given by an ideal $J$. Its radical $rad(J):=J_{red}$ presents two irreducible components: \medskip

\begin{enumerate}
	\item A component $C$, 
	with ideal
	\begin{equation} \label{JC}
	J_C=\langle b_1-4u_3,a_1-6u_3,a_2,a_3,a_4-12u_3^2,u_2 \rangle
	\end{equation}
	The ideal $J_C$ is generated by a regular
	sequence.  The affine coordinates here are $a_1,
	a_2, a_3, a_4, b_1, u_2, u_3$. 
	Therefore we have a locally complete
	intersection, which represents the curve over the
	flag $\varphi_0$ given by
	$$f=(x_0+2u_3x_1)^3, \ \ g=(x_0+2u_3x_1)^2.$$

	\item A component $E$ with ideal 
	\begin{equation} \label{JEE}
	J_E=\langle a_1,a_2,a_3,a_4,b_1 \rangle,    
	\end{equation}
	which is the whole $\mathbb{P}^2$--fiber of the
	exceptional divisor of $\mathbb{X}'$ over
	$g=x_0^2$, and $f=x_0^3$.
\end{enumerate} \medskip

Notice that there is precisely one  point in the
intersection of these two components, since
$$J_C + J_E = \langle a_1,a_2,a_3,a_4,b_1,u_2,u_3 \rangle .$$
This ideal represents the single point $(g,g',f)=(x_0^2,x_0x_1,x_0^3) \in \mathbb{Y} $. \medskip

Let's see what happens when we 
blowup
$\mathbb{Y}$ first along $C$, followed by a blowup along
$E'$, the strict transform of $E$. 

Let $\mathbb{A}^7$ be the affine neighborhood defined by $b_0=a_0=u_1=1$. The blowup of this $\mathbb{A}^7$ along $C$ is
$$\mathbb{Y}_1 = \{\left(
(\underline{a},\underline{b}),(s_0:\ldots:s_5)\right)
\ | \ s_i \cdot e_j = s_j \cdot e_i\} \subset \mathbb{A}^7 \times \mathbb{P}^5,$$
where $e_i, \ 0\leq i \leq 5$, are the equations of $J_C$, ordered as in (\ref{JC}). 

There are six choices for the local equation of the first
exceptional divisor. Choose, say,
exc$_1=u_2$, {\em i.e,} look to the affine chart $s_5=1$. 
Then, we make the following substitutions 
$$\left\{\begin{array}{rcl}
b_1-4u_3    & = & s_0u_2 \\
a_1-6u_3    & = & s_1u_2 \\
a_2         & = & s_2u_2 \\
a_3         & = & s_3u_2 \\
a_4-12u_3^2 & = & s_4u_2
\end{array}\right.
$$

The new 7 affine coordinates on the blowup are 
\begin{equation}\label{affineY1}
s_0,s_1,s_2,s_3,s_4,u_2,u_3.
\end{equation}

After the blowup $\mathbb{Y}_1 \rightarrow \mathbb{Y}$  of $\mathbb{Y}$
along $C$, 
the strict transform of the radical $J_{red}$ is given by 
\begin{equation} \label{KE'}
J_E'=\langle s_3,s_2,3s_0-2s_1,6u_3+s_1u_2,3s_4+s_1^2u_2 \rangle.
\end{equation}
It coincides with the strict transform $E'$ of $E$ under
the first blowup. Now, take $E'$ as our new blowup center and 
denote by $\mathbb{Y}_2 \rightarrow \mathbb{Y}_1$ the blowup of $\mathbb{Y}_1$ along $E'$,
$$\mathbb{Y}_2 = \{\left(
(\underline{u},\underline{s}),(t_0:\ldots:t_4)\right)
\ | \ t_i \cdot e_j = t_j \cdot e_i\} \subset \mathbb{A}^7 \times \mathbb{P}^4,$$
where $e_i, \ 0\leq i \leq 4$, are the equations of
$J_{E'}$, ordered as in (\ref{KE'}), and
$(\underline{u},\underline{s})$ as in (\ref{affineY1}).

Choose now the equation exc$_2=s_2$ in (\ref{KE'})  as
the new local exceptional equation. Equations of
$\mathbb{Y}_2$ 
become
$$\left\{
\begin{array}{rcl}
s_3           & = & t_0s_2 \\
3s_0-2s_1     & = & t_2s_2 \\
6u_3+s_1u_2   & = & t_3s_2 \\
3s_4+s_1^2u_2 & = & t_4s_2
\end{array}\right.
$$
The new 7 affine coordinates on $\mathbb{Y}_2$ are 
\begin{equation}\label{affineY2}
t_0,t_2,t_3,t_4,u_2,s_1,s_2.
\end{equation}

The two blowups so far performed are not enough to solve
the indeterminacies of our 
map $(g,f) \mapsto \omega$. However they do make    the
new scheme of indeterminacies to become reduced. 
Although reduced, the indeterminacy locus is still reducible, with two components. One of them is given by the ideal
\begin{equation} \label{JR2}
J_{R}=\langle u_2,t_3-1,t_2,t_0,3s_1-2t_4 \rangle.
\end{equation}
The ideal $J_{R}$ represents a ruled surface $R$ which is a
$\mathbb{P}^1$--subbundle of the (transform of the) exceptional divisor $\mathbb{P}(\mathcal
N_{C|\mathbb Y})$ -- note the equation exc$_1=u_2$ representing the curve $C$ and the other four linear equations on the new affine variables.

Denote by $\mathbb{Y}_3 \rightarrow \mathbb{Y}_2$ the blowup of $\mathbb{Y}_2$ along $R$. This blowup, in the affine chart $s_5=1,\ t_1=1$, is 
$$\mathbb{Y}_3 = \{\left(
(\underline{u},\underline{s},\underline{t}),(v_0:\ldots:v_4)  \right)
\ | \ v_i \cdot e_j = v_j \cdot e_i\} \subset \mathbb{A}^7 \times \mathbb{P}^4,$$
where $e_i, \ 0\leq i \leq 4$, are the equations of
$J_{R}$, ordered as in (\ref{JR2}), and
$(\underline{u},\underline{s},\underline{t})$ as in (\ref{affineY2}).

Choose now the equation exc$_3=u_2$ in (\ref{JR2})  as
the new local exceptional equation, {\em i.e,} take $v_0=1$. Equations of $\mathbb{Y}_3$ become
$$\left\{
\begin{array}{rcl}
t_3-1     & = & v_1u_2 \\
t_2       & = & v_2u_2 \\
t_0       & = & v_3u_2 \\
3s_1-2t_4 & = & v_4u_2
\end{array}\right.
$$

The new 7 affine coordinates on $\mathbb{Y}_3$ are 
\begin{equation}\label{affineY3}
v_1,v_2,v_3,v_4,u_2,s_2,t_4
\end{equation}
and the indeterminacy locus reduces to 
\begin{equation}\label{J'''L}
J_L=\langle s_2,v_1,v_2,v_3,v_4,t_4 \rangle.    
\end{equation}

The component represented by the ideal $J_L$ (\ref{J'''L}) is easy to
describe. Since the affine variables are as listed in
(\ref{affineY3}),
the six linear equations represent a
line inside the (transform of) the second blowup center $E'$ (note the presence of the equation exc$_2 = s_2$ in $J_L$), a line parametrized by the variable
$u_2$. 

We call $L$ the component described by the ideal $J_L \  
(\ref{J'''L})$. This is a reduced and irreducible 
local complete intersection. Since it is the full indeterminacy locus,
a blowup along this subscheme will solve the indeterminacies over the
neighborhood $[b_0=1, \ u_1=1],[s_5=t_1=v_0=1]$.
That is, the map to $\omega$
becomes a morphisms,
cf. \,\cite[ex. 7.17.3, p.\,168]{hartshorne2013algebraic}.

Denote by $\mathbb{Y}_4 \rightarrow \mathbb{Y}_3$ the blowup of $\mathbb{Y}_3$ along $L$,
$$\mathbb{Y}_4 = \{\left(
(\underline{u},\underline{s},\underline{t},\underline{v}),(z_0:\ldots:z_5)\right)  \ | \ z_i \cdot e_j = z_j \cdot e_i\} \subset \mathbb{A}^7 \times \mathbb{P}^5,$$
where $e_i, \ 0\leq i \leq 5$, are the equations of $J_L$, ordered as in (\ref{J'''L}), and $(\underline{u},\underline{s},\underline{t},\underline{v})$ as 
in (\ref{affineY3}).

In $\mathbb{Y}_4$ the map is solved, at least in the
affine chart [$s_5=t_1=v_0=1$].	The calculations in all
other standard neighborhoods reveal that this sequence of
four blowups, along $C$, $E$, $R$ and $L$,
will solve the map $(g,f) \mapsto \omega$ over the neighboorhod $[b_0=1, \ u_1=1]$.    

There are other five standard neighborhoods to be checked to complete an affine cover of
$\mathbb{X}'$,  corresponding to
$$[b_3=1], \ [b_2=1], \ [b_1=1], \ [b_0=1 \textrm{ and }
u_2=1], \ [b_0=1 \textrm{ and } u_3=1].$$ Similar calculations show  that these four blowups described will solve the
map over each of these neighborhoods. \medskip

The whole discussion of the indeterminacy loci was made over the fixed flag (\ref{pdv}).
Since the blowup centers as described were actually fibers of
bundles  over the variety of complete flags  ${\mathbb F}$, 
we have at the end a bundle which constitutes the
desired parameter space:
\begin{equation} \label{summarized}
\begin{array} c
\mathbb{Y}_4   \xrightarrow{\textrm{blowup } L} \mathbb{Y}_3 \xrightarrow{\textrm{blowup }R} \mathbb{Y}_2  \xrightarrow{\textrm{blowup }E} \mathbb{Y}_1  \xrightarrow{\textrm{blowup }C} \mathbb{Y} \\\\
\mathbb{Y} \xrightarrow{\mathbb{P}^4 \textrm{ bundle}} \mathbb{X}' \xrightarrow{\textrm{blowup }\mathbb{P}(\mathcal{O}_{\ensuremath{\check{\vphantom{|}\mathbb P}}^3}(-2))} \mathbb{X}
\xrightarrow{\mathbb{P}^3 \textrm{ bundle } } {\mathbb F}=\textrm{ all flags}
\end{array}
\end{equation}

This can be summarized as follows:

\subsection{\bf Theorem.}\label{Y4}
	Let $\mathbb{Y}_4$ be the variety obtained by the four
	blowups as described above. Then $\mathbb{Y}_4$ 
	is
	equipped with a morphism $\Phi$
	onto the exceptional component $E(3)$.
	\qed

\subsection{\bf Proposition.}
	The map $\Phi: \mathbb{Y}_4 \longrightarrow E(3)$ is generically injective.
	\begin{proof}\normalfont
		For a given $\omega \in E(3)$ (off the boundary) we already saw that the flag can be recovered (see Remark \ref{recupflag}). Hence, we can look at a fiber over a fixed flag. 
		
		So fix the flag $\varphi_0=(p_{0},\ell_{0},v_0) (\ref{pdv})$. The fiber of $\mathbb{X}$ at this flag is the $\mathbb{P}^3$ of special quadrics
		$$g=b_0x_0^2 + b_1x_0x_1 + b_2x_0x_2 + b_3x_1^2.$$
		Now, we will pay attention to the dense
		open set where $b_2=1$, so 
		\begin{equation} \label{gfibra}
		g=b_0x_0^2 + b_1x_0x_1 + x_0x_2 + b_3x_1^2.
		\end{equation}
		The divisibility condition (\ref{divcond1}) shows that the fiber of $\mathbb{Y}$ over a quadric as in (\ref{gfibra}) is the $\mathbb{P}^4=\{(a_0:a_1:a_2:a_3:a_5) \}$ of special cubics
		$$f=(a_0x_0+a_1x_1 +a_2x_2+a_3x_3)x_0^2
		+ a_5\left(b_1x_0x_1^2 +    x_0x_1x_2 +
		\frac{2}{3}   b_3x_1^3\right)
		.$$
		Computing $\omega = \dfrac{3fdg-2gdf}{x_0}$ one can see that there are no indeterminacies on this neighborhood $[b_2=1]$. Moreover,
		$$2g\dfrac{\partial f}{\partial x_3} - 3f\dfrac{\partial g}{\partial x_3} = 2a_3gx_0^2\ ,$$ 
		and the quadric $g$ can be recovered from
		$\omega $ just by looking at the coefficient of $dx_3$ (at least on the dense open set $a_3=1$).
		
		The coefficient of $dx_2$ in $\omega$ is
		
		\begin{eqnarray*} 
			\dfrac{1}{x_0}\left(2g\dfrac{\partial f}{\partial x_2} -
			3f\dfrac{\partial g}{\partial x_2}\right) = 
			(2a_2b_0-3a_0)x_0^3+(2a_5b_0+2a_2b_1-3a_1)x_0^2x_1+ \\ + (-a_5b_1+2a_2b_3)x_0x_1^2
			-a_2x_0^2x_2-a_5x_0x_1x_2-3a_3x_0^2x_3.
		\end{eqnarray*}
		
		Hence the coefficients $a_2$ and $a_5$ can be recovered from $x_0^2x_2dx_2$ and $x_0x_1x_2dx_2$, 
		respectively. And since the coefficients $b_0$ and $b_1$ are also known we obtain the coefficients $a_0$ and $a_1$ from $x_0^3dx_2$ and $x_0^2x_1dx_2$ 
		respectively.
	\end{proof} 

\section{Calculation of Degree}\label{deg}

By construction, $\mathbb{Y}_4$ is equipped with a line bundle $\mathcal{W}$, pullback of the line bundle $\mathcal{O}_{\mathbb{P}^{44}}(-1)$ over 
${\mathbb{P}^{44}}=\mathbb{P}(H^0(\Omega^1_{\mathbb{P}^3}(4)))$. The fiber of $\mathcal{W}$ over a point in $\mathbb{Y}_4$ is the rank one space spanned by the computed 1-form $\omega$.

$$\begin{array}{ccccc}
\mathcal{W}=\Phi^*(\mathcal{O}(-1))&   & \mathcal{O}(-1)  & &  \\
\downarrow &  & \downarrow  & & \\
\mathbb{Y}_4 & \stackrel{\Phi}\longrightarrow    & \mathbb{P}(H^0(\Omega^1_{\mathbb{P}^3}(4))) & \supset & E(3)=\Phi(\mathbb{Y}_4)\\
\end{array}$$

\subsection{\bf Theorem.} \label{integral}
	The degree of the exceptional component of
	codimension one and degree two foliations
	in $\mathbb{P}^3$ is given by
	$$\int\limits_{\mathbb{Y}_4}-c_1^{13}(\mathcal{W}) \cap [\mathbb{Y}_4].$$
	\begin{proof}\normalfont
		We have $\dim (\mathbb{Y}_4) =13$. Since the map $\Phi$ is generically injective the required degree is (cf. definition of deg$_f\tilde{X}$ in \cite{Fulton}, page 83) $$\int\limits_{\mathbb{Y}_4}c_1(\Phi^*\mathcal{O}(1))^{13} \cap [\mathbb{Y}_4] = \int\limits_{\mathbb{Y}_4}-c_1(\Phi^* \mathcal{O}(-1))^{13} \cap [\mathbb{Y}_4]= \int\limits_{\mathbb{Y}_4}-c_1(\mathcal{W})^{13} \cap [\mathbb{Y}_4].$$
	\end{proof}

In order to obtain the value of the integral in Theorem \ref{integral}, we apply Bott's formula
\begin{equation}  \label{BOTT}
\int\limits_{\mathbb{Y}_4}-c_1^{13}(\mathcal{W}) \cap [\mathbb{Y}_4]=\sum\limits_{F}\dfrac{-c_1^T(\mathcal{W})^{13} \cap [F]_T}{c_{top}^T(\mathcal{N}_{F|\mathbb{Y}_4})},\end{equation}
where the sum runs through all fixed components $F$ under
a convenient action of the torus $\mathbb{C}^*$.  The $\mathcal{N}_{F|\mathbb{Y}_4}$ appearing in the denominator denotes
the normal bundle of a fixed component $F$ in $\mathbb Y_4$.
Fix  an action 
\begin{equation}\label{actionp3}
\begin{array}{ccc}
\mathbb{C}^* \times \ensuremath{\check{\vphantom{|}\mathbb P}}^3 & \longrightarrow & \ensuremath{\check{\vphantom{|}\mathbb P}}^3 \\
(t,x_i) & \mapsto & t^{w_i}x_i
\end{array}
\end{equation}
of $\mathbb{C}^*$ in $\ensuremath{\check{\vphantom{|}\mathbb P}}^3$ with distinct weights $w_i$, $i \in \{0,1,2,3 \}$.
The only fixed flags are the 24 standard ones
$$p_{ijk}=\{x_i=x_j=x_k=0\} \in \ell_{ij}=\{x_i=x_j=0\} \subset v_i=\{x_i=0\}.$$

Over each one of these fixed flags, we can find 72 fixed isolated points and 5 fixed lines. So, there is a total of $72 \cdot 24
= 1728$ fixed points and $5 \cdot 24 = 120$ fixed lines. 

A detailed exposition of the fixed points and the computations of their contributions on Bott's formula (\ref{BOTT}), including scripts for  \textit{Macaulay2} for this sum and for the resolution of singularities also in the other neighborhoods can be found in  

\url{https://sites.google.com/a/ifsudestemg.edu.br/nucleo-de-matematica/artur/trabalhos/Exceptional.pdf}  

\subsection{\bf Theorem.} \label{GRAUE3}
	The degree of the exceptional component of foliations of codimension one and degree two in $\mathbb{P}^3$ is \textbf{168208}. \medskip
	
	In the next Table \ref{tab} we list the components of
	the space of foliations of degree two and codimension one
	in $\mathbb{P}^3$ as described in  \cite{Alcides}, and
	their respective degrees.
	Notice that the degree of the logarithmic
	components have not been found yet in the literature to
	the best of our knowledge.

\begin{table}
	\caption{Components of $\overline{\mathcal{F}(2,3)}$ and their Degrees}
	\label{tab}       
	\begin{tabular}{lll}
		\hline\noalign{\smallskip}
		Component & Degree & Reference  \\
		\noalign{\smallskip}\hline\noalign{\smallskip}
		Linear Pullbacks $\mathcal{S}(2,3)$ & 1320 & Ferrer and Vainsencher (to appear)  \\
		Rational $\mathcal R(2,2)$ & 1430 & \cite{CukiermanVI} \\
		Rational $\mathcal R(1,3)$ & 700 & \cite{CukiermanVI} \\
		Logarithmic $\mathcal{L}(1,1,1,1)$ & ???? &  \\
		Logarithmic $\mathcal{L}(1,1,2)$ & ???? &  \\
		Exceptional Component $E(3)$ & 168208 & Theorem \ref{GRAUE3} \\
		\noalign{\smallskip}\hline
	\end{tabular}
\end{table}

\medskip

\section{A geometric interpretation of the degree}\label{geo}

For a codimension one foliation in $\mathbb{P}^n$, given by the differential form
$$\omega = \sum_{i=0}^{n} A_idx_i,$$
the hyperplane defined by the distribution at a point $p \in \mathbb{P}^n$ is 
$$H=\sum_{i=0}^{n} A_i(p)x_i = 0.$$

Thus, a tangent direction $v=(v_0: \ldots : v_n) \in \mathbb{P}(T_p\mathbb{P}^3)$ lies on this hyperplane $H$ if 
$$\sum_{i=0}^{n} A_i(p) \cdot v_i = 0.$$
This can be thought of as a linear equation on the
coefficients of the $A_i$.
Hence, the point $(p,v) \in \mathbb{P}(T\mathbb{P}^3) $ 
defines a hyperplane in the projective space of distributions.
\begin{equation}\label{condlinear}
\omega(p) \cdot v = 0.
\end{equation}

The equation (\ref{condlinear}) shows that the degree of
a $m$--dimensional component of the space of codimension
one foliations in $\mathbb{P}^n$ can be interpreted as
the number of such foliations that are tangent  to $m$ 
general
directions in $\mathbb{P}^n$.  

In particular, Theorem \ref{GRAUE3} means that there are 168208 exceptional foliations tangent to 13 general directions in $\mathbb{P}^3$.

\subsection*{\bf Acknowledgements}
	Our hearts and minds go with
	Flaviano Bahia P. Vieira (1984-2011). He has worked on this subject,
	especially  during a visit at the UVa. His untimely
	death has  deprived us of the company of a beloved friend.

\medskip
\hfill$\begin{array}l
\text{IFSudesteMG -- JF, MG -- Brasil}\\
\text{artur.rossini@ifsudestemg.edu.br}
\end{array}$ \hfill
$\begin{array}l
\text{UFMG  --  BH, MG - Brasil}\\
\text{israel@mat.ufmg.br}
\end{array}$

\end{document}